\documentclass[12pt]{amsart}
\usepackage{amscd, amsfonts, amssymb}
\usepackage{latexsym}
\usepackage{verbatim}
\usepackage{enumerate}

\newcommand\ra{\rightarrow}

\newcommand{\RR}{{\mathbb R}}

\DeclareMathOperator{\Aut}{Aut}

\DeclareMathOperator{\GL}{GL}
\DeclareMathOperator{\Gal}{Gal}

\numberwithin{equation}{section}

\newtheorem{thm}[equation]{Theorem}

\theoremstyle{definition}
\newtheorem{defn}[equation]{Definition}

\theoremstyle{remark}
\newtheorem{rem}[equation]{Remark}
\theoremstyle{remark}
\newtheorem{rems}[equation]{Remarks}

\newcommand{\ovl}{\overline}


\title[Complete reducibility and separable field extensions]
{Complete reducibility and \\ separable field extensions}

\author[M.\  Bate]{Michael Bate}
\address
{Department of Mathematics,
University of York,
York YO10 5DD,
United Kingdom}
\email{meb505@york.ac.uk}

\author[B.\ Martin]{Benjamin Martin}
\address
{Mathematics and Statistics Department,
University of Canterbury,
Private Bag 4800,
Christchurch 8140,
New Zealand}
\email{B.Martin@math.canterbury.ac.nz}

\author[G. R\"ohrle]{Gerhard R\"ohrle}
\address
{Fakult\"at f\"ur Mathematik,
Ruhr-Universit\"at Bochum,
44780 Bochum, Germany}
\email{gerhard.roehrle@rub.de}

\begin{document}

\begin{abstract}
Let $G$ be a connected reductive linear algebraic group.
The aim of this note is to settle a question of
J-P.\ Serre concerning the behaviour of his notion of
$G$-complete reducibility under
separable field extensions.
Part of our proof relies on the recently
established Tits Centre Conjecture
for the spherical building of
the reductive group $G$.
\end{abstract}

\maketitle

\section{Introduction}
\label{sec:intro}

Throughout, $G$ denotes a connected reductive
linear algebraic group defined over a field $k$.
Following Serre, \cite{serre2}, a
subgroup $H$ of $G$
is called \emph{$G$-completely reducible over $k$ ($G$-cr over $k$)}
if whenever $H$
is contained in a $k$-defined parabolic subgroup $P$ of $G$, there exists
a $k$-defined
Levi subgroup of $P$ containing $H$.
In case $V$ is a finite dimensional $k$-vector space and $G = \GL(V)$,
a subgroup $H$ of $G$ is $G$-completely reducible over $k$
precisely when $V$ is a semisimple $H$-module, \cite[1.3, 3.2.2]{serre2}.
In this sense, Serre's notion generalizes the usual
concept of complete reducibility in representation theory.
For more details and further
results on this notion,
see \cite{serre1.5}, \cite{serre2}, \cite{BMR}, \cite{BMRT}, and \cite{GIT}.

The following theorem answers a question of Serre.

\begin{thm}
\label{thm:serrequestion}
Suppose $k_1/k$ is a separable extension of fields.
Let $G$ be a reductive group defined over $k$, and let $H$ be
a $k$-defined subgroup of $G$.
Then $H$ is $G$-completely reducible over $k$ if and only if $H$ is
$G$-completely reducible over $k_1$.
\end{thm}

The reverse implication in Theorem \ref{thm:serrequestion}
is proved in \cite[Thm.\ 5.11]{GIT}.
The proof of \cite[Thm.\ 5.11]{GIT}
rests on a general rationality result, \cite[Thm.\ 3.1]{GIT}, concerning
$G$-orbits in an affine variety.
We present a proof of the forward direction of the statement
in Section \ref{sec:proof} based on the recently established
Tits Centre Conjecture, Theorem \ref{conj:tcc}.

\begin{rems}
\label{rems:thm}
(i).
In \cite[Ex.\ 5.12]{GIT}, we showed that
Theorem \ref{thm:serrequestion} holds when $G = \GL(V)$.

(ii).
Theorem \ref{thm:serrequestion}
was proved in \cite[Thm.~5.8]{BMR} for
$k$ perfect, by passing back and forth between $k$ and its algebraic closure
$\ovl{k}$ and between $k_1$ and $\ovl{k}$.
In general this approach fails,
because the extension $\ovl{k}/k$ need not be separable.

(iii).
There are examples showing that each implication in
Theorem \ref{thm:serrequestion}
fails without the separability assumption on the extension
$k_1/k$; see \cite[Ex.~5.11]{BMR} and \cite[Ex.~7.22]{BMRT}.
%
\end{rems}

\section{The Centre Conjecture for spherical buildings}
\label{sec:tcc}

Let $\Delta_k$ denote the spherical building of $G$ over $k$, \cite[Sec.\ 5]{tits1}:
the simplices of $\Delta_k$ correspond to $k$-defined parabolic subgroups
of $G$. 
Given a $k$-defined parabolic subgroup $P$ of $G$, we denote
the simplex corresponding to $P$ in $\Delta_k$ by $\sigma_P$.
Throughout, we identify $\Delta_k$ with its geometric realization,
which is a bouquet of spheres \cite{serre2}.

An \emph{apartment} in $\Delta_k$ consists of the simplices $\sigma_P$
corresponding to
all $k$-defined parabolic subgroups $P$ of $G$
which contain a fixed maximal $k$-split torus of $G$;
it is a subcomplex whose geometric realization is a sphere.
Any two points of $\Delta_k$ lie in a common apartment.
We say that $x,y\in \Delta_k$ are
\emph{opposite} if they are opposite
in some apartment that contains them both.
It can be shown that if $x$ and $y$ are opposite
in some apartment that contains them both, then
they are opposite in any apartment that contains them both.
If $x,y\in \Delta_k$ are not opposite,
then there is a unique geodesic joining them, \cite[\S 2.1.4]{serre2}.
Two simplices $\sigma_P$ and $\sigma_Q$ are said to be \emph{opposite}
if every point of  $\sigma_P$ is opposite a point of $\sigma_Q$, \cite[\S 2.1.4]{serre2}.
In terms of parabolic subgroups of $G$,
the simplices $\sigma_P$ and $\sigma_Q$ corresponding to $k$-defined
parabolic subgroups $P$ and $Q$ of $G$
are opposite in $\Delta_k$ if and only
if $P\cap Q$ is a common Levi subgroup of $P$ and $Q$
(this Levi subgroup is then automatically $k$-defined).


Let $\Sigma$ be a subcomplex of $\Delta_k$.
We say that $\Sigma$ is \emph{convex}
if whenever $x,y\in \Sigma$ are not opposite, then $\Sigma$
contains the geodesic between $x$ and $y$, \cite[\S 2.1]{serre2}.

Suppose $\Sigma$ is a convex subcomplex of $\Delta_k$.
Serre has shown that $\Sigma$ is \emph{contractible} --- that is, $\Sigma$
has the homotopy type of a point --- if and only if there exists a point
of $\Sigma$ which has no opposite in $\Sigma$; see \cite[\S 2.2]{serre2}.
The following terminology is due to Serre \cite[Def.\ 2.2.1]{serre2}:

\begin{defn}
\label{def:Delta-cr}
Let $\Sigma$ be a convex subcomplex of $\Delta_k$.
We say that $\Sigma$ is \emph{$\Delta_k$-completely reducible}
(or $\Delta_k$-cr) if every simplex in $\Sigma$
has an opposite in $\Sigma$.
\end{defn}

Serre has shown that the group-theoretic definition of
$G$-complete reducibility over $k$ has the following
building-theoretic interpretation,
\cite{serre2}:
Given a subgroup $H$ of $G$, let
$$
\Delta_k^H = \{\sigma_P \mid P
\textrm{ is a $k$-defined parabolic subgroup containing } H\}.
$$
Then $\Delta_k^H$ is a convex subcomplex of $\Delta_k$ (\cite[Prop.\ 3.1]{serre2}),
the \emph{fixed point subcomplex} of $\Delta_k$ under the action of $H$,
and $H$ is $G$-completely reducible over $k$ if and only if
$\Delta^H_k$ is $\Delta_k$-cr,
\cite[2.3.1, 3.2]{serre2}.
Equivalently, $H$ is \emph{not} $G$-completely reducible over $k$
if and only if $\Delta_k^H$ is
contractible.

\begin{defn}
\label{defn:centre}
Let $\Sigma$ be a subcomplex of $\Delta_k$ and let $x\in \Sigma$.
Let $\Gamma$ be a group which acts on $\Delta_k$ by means of building
automorphisms, \cite{tits1}, i.e., suppose there is a homomorphism
$\Gamma \to \Aut \Delta_k$, where $\Aut \Delta_k$ is the group of building automorphisms of $\Delta_k$.
We say that $x$ is a \emph{$\Gamma$-centre of $\Sigma$}
if $x$ is fixed by any element of $\Gamma$ that stabilizes $\Sigma$ setwise.
\end{defn}

The following theorem is known as the ``Centre Conjecture'' of J.\ Tits,
cf.\  \cite[Lem.\ 1.2]{tits0}, \cite[\S 4]{serre1.5},
\cite[\S 2.4]{serre2}, \cite{tits2}, \cite[Ch.\ 2, \S 3]{mumford},
\cite[Conj.\ 3.3]{rousseau}.
It has recently been proved
in a series of intricate case-by-case arguments
by B.\ M\"uhlherr and J.\ Tits \cite{muhlherrtits}
($G$ of classical type or type $G_2$),
B.\ Leeb and C.\ Ramos-Cuevas \cite{lrc}
($G$ of type $F_4$ or $E_6$)
and C.\ Ramos-Cuevas \cite{rc}
($G$ of type $E_7$ or $E_8$).

\begin{thm}[Tits' Centre Conjecture]
\label{conj:tcc}
Let $\Sigma$ be a convex contractible subcomplex of $\Delta_k$.
Then $\Sigma$ has an $\Aut \Delta_k$-centre.
\end{thm}

\begin{rem}
\label{rem:tcc}
Suppose $G$ is semisimple 
and $k$ is a perfect field.
It follows from  \cite[5.7.2]{tits1} that $\Aut G$
is an algebraic group also defined over $k$.
In \cite[Thm.\ 5.31]{GIT}, we give a uniform proof
of the following special case of the Centre Conjecture:
Let $H$ be a subgroup of $G$.
If $\Delta_k^H$ is contractible, i.e., if $\Delta_k^H$ is not
$\Delta_k$-cr, then $\Delta_k^H$ admits an $(\Aut G)(k)$-centre.
The proof of this result in \cite{GIT} utilizes methods from
geometric invariant theory and the concept of optimal destabilizing
parabolic subgroups.

Let $k$ be a field, let
$k_s$ denote its separable closure,
and let $\ovl{k}$ denote its algebraic closure.
Note that $k_s = \ovl{k}$ if $k$ is perfect.
Thanks to \cite[5.7.2]{tits1},  $\Gamma : = \Gal(k_s/k)$ acts on
$\Delta_{\ovl{k}}$ via building automorphisms.
In \cite[Thm.\ 5.33]{GIT}, we show that if 
$H$ is a $k$-defined subgroup of $G$ such that
$\Delta_{\ovl{k}}^H$ is contractible,  then $\Delta_{\ovl{k}}^H$ admits a
$\Gamma$-centre.
The proof of \cite[Thm.\ 5.33]{GIT} rests on a
rationality result concerning $G$-cr subgroups of $G$,
\cite[Prop.\ 5.14(iii)]{GIT}.

Both \cite[Thm.\ 5.31]{GIT} and \cite[Thm.\ 5.33]{GIT}
improve on \cite[Thm.\ 3.1]{BMR:tits}.
\end{rem}

\section{Proof of Theorem \ref{thm:serrequestion}}
\label{sec:proof}

As noted above the reverse implication of
Theorem \ref{thm:serrequestion}
is proved in \cite[Thm.\ 5.11]{GIT}.
We deduce the other direction with the aid of
Theorem \ref{conj:tcc}.

Suppose $k_1/k$ is an algebraic separable extension of fields and
let $\Delta_{k}$ and $\Delta_{k_1}$ denote the buildings of
$G$ over $k$ and $k_1$, respectively.
By the reverse implication, one may suppose that $k_1/k$ is Galois.
Then the Galois group $\Gamma := \Gal(k_1/k)$ acts simplicially
on $\Delta_{k_1}$, i.e., $\Gamma$ permutes the
set of $k_1$-defined parabolic subgroups of $G$.
Moreover, the subcomplex of $\Delta_{k_1}$ consisting of $\Gamma$-stable
simplices is just $\Delta_k$.

It is convenient to reduce to the case when $H$ is not
contained in any $k$-defined
Levi subgroup of any proper $k$-defined parabolic subgroup of $G$.
To do this, 
we let $L$ be minimal such that $L$ is a $k$-defined Levi subgroup of some 
$k$-defined parabolic subgroup $P$ of $G$ and $H\subseteq L$. 
Then $L$ is also $k_1$-defined, and by a result of
Serre, \cite[Prop.\ 3.2]{serre2},
$H$ is $G$-completely reducible over $k$ (resp.\ $k_1$)
if and only if $H$ is $L$-completely reducible over $k$ (resp.\ $k_1$).
Now if $L'$ is a $k$-defined Levi subgroup of some
proper $k$-defined parabolic subgroup $Q$ of $L$,
then $QR_u(P)$ is a $k$-defined parabolic
subgroup of $G$, and $L'$ is a Levi subgroup of $QR_u(P)$,
\cite[Prop.\ 4.4]{boreltits}.
Since $L$ is minimal among those $k$-defined Levi subgroups 
of $k$-defined parabolic subgroups of $G$ that contain $H$, 
$H$ cannot be contained in $L'$. 
By replacing $G$ with $L$, we can now assume that $H$ is not contained
in any $k$-defined Levi subgroup of any proper $k$-defined parabolic subgroup of $G$.

Suppose that $H$ is not $G$-completely reducible over $k_1$.
Then $\Delta_{k_1}^H$ is contractible, and since $H$ is $k$-defined,
$\Delta_{k_1}^H$ is $\Gamma$-stable.
Since $\Delta_{k_1}^H$ is a convex contractible subcomplex of $\Delta_{k_1}$,
it follows from Theorem \ref{conj:tcc}
that $\Gamma$ fixes a point
of $\Delta_{k_1}^H$, and this point lies in some (minimal) simplex $\sigma_P$,
where $P$ is a proper $k_1$-defined parabolic subgroup of $G$.
Since the action of $\Gamma$ on $\Delta_{k_1}$ is simplicial,
$P$ is stabilized by $\Gamma$,
which is equivalent to saying that $P$ is $k$-defined.
Now, by assumption, $H$
is not contained in any $k$-defined Levi subgroup of $P$,
so $H$ is not $G$-completely reducible over $k$.
This completes the proof of Theorem \ref{thm:serrequestion}.

\begin{rem}
\label{rem}
In  \cite[Thm.\ 4.13]{BMRT:relative}, we prove a generalization of
the reverse implication of Theorem \ref{thm:serrequestion}
in the setting of ``relative complete reducibility''.
The arguments above 
used to derive the
forward direction of Theorem \ref{thm:serrequestion}
do not apply to this more general situation, as the
relevant subset 
in $\Delta_{k_1}$
is only a convex subset but not a subcomplex of $\Delta_{k_1}$.
Thus Theorem \ref{conj:tcc} does not apply.
%
\end{rem}


\bigskip
{\bf Acknowledgements}:
The paper was written during a stay of the first and third authors
at the Max Planck Institute for Mathematics in Bonn.
The authors also acknowledge financial support
from the DFG-priority program SPP1388 ``Representation Theory'', 
Marsden Grant UOC0501,
and The Royal Society.

\bigskip

\end{document}